**A hysteresis model for two-dimensional input signals**

S. A. Belbas[*] and Young Hee Kim[#]

\* Mathematics Department, University of Alabama, Tuscaloosa, AL. 35487-0350. USA.
\# Department of Mathematics, Tuskegee University, Tuskegee, AL. 36088. USA.

Corresponding author: S. A. Belbas.
e-mail: SBELBAS@BAMA.UA.EDU

Abstract.   We formulate and analyze a new model of vector hysteresis for the case of two-dimensional input signals. We prove certain properties of this model and we present the solutions to two identification problems connected with our model.

1. Introduction

   The phenomenon of hysteresis appears in many different kinds of applications: the theory of magnetism is the prototypical instance of hysteresis [M1, BM]. However, hysteresis also appears in many other areas, for example in the mathematical modeling of viscoelastic materials, in economics [C, DCP, PCGL, TK, P], in biology, in the modern topic of sociophysics and econophysics, in hydrology [FMKP], etc.

   The best definition of hysteresis that we know of is due to R. Cross: hysteresis is a phenomenon that exhibits <u>selective memory</u> and <u>remanence</u>. <u>Selective memory</u> means that the output of a hysteresis operator depends only on certain particular characteristics of the past history of the input, and not on the complete past history; <u>remanence</u> means that a cyclic input, starting from one value and eventually returning to the same value, results in a nonzero net effect in the hysteretic output. Ref. [P] introduces another requirement for hysteresis: the memory should be erasable, i.e. certain conditions on the behavior of the input signal should erase all or part of the previous history of the hysteretic output. The model we introduce in this paper also has the erasability property, which, following [M1], we have termed "wiping-out property".

   In the theory of hysteresis for one-dimensional input signals and one-dimensional output signals, a prominent role is played by the Preisach model [M1, M2, KrP, V]. The model of Preisach was originally developed in the context of the theory of magnetism, but later (starting with the work of Krasnosel'skii and Pokrovskii) was analyzed as a purely mathematical model.

   The research literature contains many works that extend the Preisach model to the case of a vector input. Most of these models are generally motivated by specific applications, usually to the theory of magnetism [M3]. Our model is different from other models of vector hysteresis.

   In [DCP], a Preisach model (with scalar input) has been introduced to model the evolution of the amount of active firms in an economy as a hysteretic function of economic growth; it is plausible that the amount of active firms can in general depend hysteretically on more than one factors, for example not only on the economic growth but also on the amount of imports. In magnetostriction, strain depends hysteretically on two quantities, stress and magnetic field (see [M1]). In hydrology, the pore content of water in porous media has been modelled as a hysteretic function of piezometric head ([FMKP] and the further references therein); other sources, for example [KP], report hysteretic relations among more than two quantities, although Preisach-type hysteresis is not considered in [KP], and consequently hysteresis models with vector-valued inputs are of interest in hydrology.

   In this paper, we have been interested in developing a general mathematical model of hysteresis for two-dimensional input signals that extends in a natural way the classic Preisach model. Our model is based on the following basic idea: the pairs of scalar parameters that define a family of non-ideal relays in the scalar Preisach model are replaced by a collection of pairs of curves in our two-dimensional vector model. Because of the high level of degrees of freedom in choosing pairs of curves, we obtain a model that covers a great variety of situations and thus, for cases to which it can be applied, it has greater flexibility to accommodate experimental data. The model introduced in this paper poses new mathematical questions which have no counterpart in the theory of standard Preisach hysteresis. Among these are: the canonical representation of our hysteresis operator via two signals that obey certain compatibility restrictions, and the identification problem which concerns not only the weight function but also the two families of curves that determine the two-dimensional analogue of non-ideal relays.

2. Definition and simple properties of a two-dimensional hysteresis operator

We start by defining the concept of a non-ideal relay when the input signal takes values in $\mathbb{R}^2$. It will be seen that many ideas extend also to the case of input signals taking values in $\mathbb{R}^n$, $n \geq 3$.

Let $\Omega$ be an open subset of $\mathbb{R}^2$. The class $I_2(\Omega)$ of input signals consists of continuous functions $u$ from a time interval $[0,T]$ into $\Omega$; we shall write $u(t) = \begin{bmatrix} u_1(t) \\ u_2(t) \end{bmatrix}$ for such signals.

Definition 2.1. Let $D_0$, $D_1$ be two open subsets of $\Omega$ such that $D_0 \cup D_1 = \Omega$. We denote by $(\gamma_0)$, $(\gamma_1)$ the relative boundaries of $D_0$, $D_1$, respectively, i.e.

$$(\gamma_i) := \Omega \cap \partial D_i, \ i = 0, 1$$

(2.1)

We assume that $(\gamma_0)$, $(\gamma_1)$ are continuous curves in $\mathbb{R}^2$.

Let $u$ be a signal of the class $I_2(\Omega)$. A relay operator $R$ is defined as follows:
the possible values of $R$ are *0* and *1*;
at the initial time $t = 0$, $(Ru)(0) = i \Rightarrow u(0) \in D_i$;
the function $(Ru)(t)$, as a function of *t*, is continuous from the right and has limits from the left for every $t \in [0,T]$; the value of $(Ru)(t)$ changes at the moments of exit of $u$ from the sets $D_i$; if $(Ru)(t^-) = i$ and $u(t) \in (\gamma_i)$, then $(Ru)(t) = (i+1)_{\mod 2}$. ///

For any open set $D \subseteq \Omega$ and any $t \in [0,T)$, we define the underline{exit time of *u* from *D* at times subsequent to *t*} by

$$\tau(t, u, D) := \inf\{s : s > t \ \& \ u(s) \notin D\}$$

(2.2)

According to definition 2.1, the values of the relay $(Ru)(t)$ will switch (from *0* to *1* or from *1* to *0*) at each exit moment that is subsequent to the previous exit moment.
We have:

Theorem 2.1. If the minimum distance between the curves $(\gamma_0)$ and $(\gamma_1)$ remains bounded away from *0*, i.e.

$$\inf\{|x-y|: x \in (\gamma_0), y \in (\gamma_1)\} \equiv \delta > 0 \tag{2.3}$$

then there is a finite collection $\tau_1 < \tau_2 < \cdots < \tau_N$ of exit times such that

$$\tau_1 = \tau(0, u, D_{i_0}), \ i_0 \equiv (Ru)(0); \ \tau_j = \tau(\tau_{j-1}, u, D_{(i_0+j-1) \bmod 2}) \tag{2.4}$$

Proof. Suppose $\tau_j$ is an exit moment as indicated in (2.2), and purely for notational simplicity, assume that $\tau_j$ is a moment of exit of $u$ from $D_0$.
Since $u$ is continuous, it is uniformly continuous on $[0,T]$, and consequently, for every $\varepsilon > 0$, there exists an $\eta = \eta(\varepsilon)$ such that $|t'-t''| < \eta(\varepsilon) \Rightarrow |u(t')-u(t'')| < \varepsilon$.
Let $\delta$ be defined as in (2.3).
Then $\tau(\tau_j, u, D_1) \geq \tau_j + \eta(\delta)$; we set $\tau_{j+1} \equiv \tau(\tau_j, u, D_1)$.
We have just shown that $\tau_{j+1} - \tau_j \geq \eta(\delta)$, thus there is a finite collection (possibly empty) of exit moments of $u$ from the sets $D_i$ in the time-interval $[0,T)$. ///

Corollary 2.1. If $u$ is Lipschitz with Lipschitz constant $M$, that is,

$$|u(t')-u(t'')| \leq M|t'-t''| \ \text{ for all } t', t'' \text{ in } [0,T], \text{ then } \tau_{j+1} - \tau_j \geq \frac{\delta}{M} . \ ///$$

In order to build an analogue of the classical Preisach operator, we consider two families of open subsets of $\Omega$, say $D_0(c_0), D_1(c_1)$, each family being a one-parameter family, and the corresponding curves $(\gamma_0(c_0)), (\gamma_1(c_1))$.(See Figure 2.1)

Definition 2.2. A pair of parameters $(c_0, c_1)$ will be called <u>admissible</u> if

$$D_0(c_0) \cup D_1(c_1) = \Omega .$$

Each admissible pair $(c_0, c_1)$ and initial value $i_0$ that satisfies the conditions of definition 2.1, defines a relay operator relative to the sets $D_0(c_0), D_1(c_1)$; we shall denote that relay operator by $R^{c_0 c_1}$.
Let $C \subseteq IR^2$ be the set of all admissible pairs $(c_0, c_1)$.
A collection of initial values for the relays $R^{c_0 c_1}$, say $i_0(c_0, c_1), (c_0, c_1) \in C$, will be called <u>admissible</u> if the function $i_0 : C \to \{0, 1\}$ is Lebesgue measurable. ///

Definition 2.3. For each families of curves $(\gamma_i(c_i))$, the underline{partial order} $\underset{(i)}{\prec}$ is defined by

$$(\gamma_i(c_i)) \underset{(i)}{\prec} (\gamma_i(c_i')) \quad \text{if } D_i(c_i) \subseteq D_i(c_i'), \quad i = 0, 1. ///$$

In order to further analyze the relay operators and obtain a useful representation of the relays, we assume the following properties:

(1) Every point $x \in \Omega$ belongs to exactly one curve of the type $(\gamma_0)$, and also to exactly one curve of the type $(\gamma_1)$. We set $c_0(x)$ for the parameter $c_0$ such that $x \in (\gamma_0(c_0(x)))$ and $c_1(x)$ for the parameter $c_1$ such that $x \in (\gamma_1(c_1(x)))$.

(2) Each of the partial orders $\underset{(i)}{\prec}$ is a total order, that is, for every two sets $D_i(c_i)$ and $D_i(c_i')$, with $c_i \neq c_i'$, we have exactly one of the two conditions:

$$D_i(c_i) \underset{\neq}{\subseteq} D_i(c_i'), \text{ or } D_i(c_i') \underset{\neq}{\subseteq} D_i(c_i).$$

(3) The parameters $c_0$ and $c_1$ are ordered according to the total order of condition (2), that is,

$$c_i < c_i' \quad \text{if } D_i(c_i) \underset{\neq}{\subseteq} D_i(c_i') \quad (i = 0, 1).$$

(4) For every two points $x$, $x'$ in $\Omega$,

$$c_0(x) < c_0(x') \Leftrightarrow c_1(x) > c_1(x').$$

Under conditions (1) - (4) above, we shall prove the following theorem:

Theorem 2.2. Given an input signal $u(t)$, we define two scalar-valued signals $K_0(t)$ and $K_1(t)$ by $K_i(t) = c_i(u(t))$ ($i = 0, 1$), where $c_i(\cdot)$ is as in condition (1) above. (See Figure 2.2)
Then the relay operators $R^{c_0 c_1}$ can be expressed as follows:

(i) If $K_0(t) < c_0$ and $K_1(t) < c_1$, then $(R^{c_0 c_1} u)(t) = (R^{c_0 c_1} u)(t^-)$;
(ii) It is not possible to have both $K_0(t) > c_0$ and $K_1(t) > c_1$;
(iii) If $K_0(t) > c_0$ and $K_1(t) < c_1$, then $(R^{c_0 c_1} u)(t) = 1$;

(iv) If $K_0(t) < c_0$ and $K_1(t) > c_1$, then $(R^{c_0 c_1} u)(t) = 0$.

Proof. (i) The condition $K_0(t) < c_0$ implies that $u(t) \in D_0(c_0)$; similarly, the condition $K_1(t) < c_1$ implies $u(t) \in D_1(c_1)$. Thus, when we have both $K_0(t) < c_0$ and $K_1(t) < c_1$, we also have $u(t) \in D_0(c_0) \cap D_1(c_1)$, and therefore $(R^{c_0 c_1} u)(t)$ is continuous at $t$.

(ii) If we have $K_0(t) > c_0$, then $u(t) \notin D_0(c_0)$; similarly, if we have $K_1(t) > c_1$, then $u(t) \notin D_1(c_1)$.
So if we have both $K_0(t) > c_0$ and $K_1(t) > c_1$, then $u(t) \notin D_0(c_0) \cup D_1(c_1)$.
This contradicts the condition of admissibility of the pair $(c_0, c_1)$, since $u(t) \in \Omega$.

(iii) If $K_0(t) > c_0$, then $u(t) \notin D_0(c_0)$; if also $K_1(t) < c_1$, then $u(t) \in D_1(c_1)$.
Therefore, $u(t) \in D_1(c_1) \setminus D_0(c_0)$, and consequently $(R^{c_0 c_1} u)(t) = 1$.

(iv) The proof of (iv) is entirely analogous to the proof of (iii). ///

Corollary 2.2. Under conditions (1)-(4), the relay operator $(R^{c_0 c_1} u)(t)$, for $t \in [t_0, T]$, is completely determined by the value of the relay at $t = t_0$ and the history of the signal $K(t)$.

Proof. From cases (i), (iii), and (iv) of Theorem 2.2, the value of $(R^{c_0 c_1} u)(t)$ is determined on the basis of conditions that are described by using the signal $K(t)$.
In order to complete all possibilities for the signal $K(t)$,
we have $(R^{c_0 c_1} u)(t) = 1$ if $K_0(t) = c_0$ and $K_1(t) < c_1$, since the relay $R^{c_0 c_1}$ satisfies, by Definition 2.1, $(R^{c_0 c_1} u)(t) = 1$ when $u(t) \in \gamma_0(c_0)$.
Also, by conditions (1)-(4), the condition $K_0(t) = c_0$ is equivalent to $u(t) \in \gamma_0(c_0)$.
Similarly, we have $(R^{c_0 c_1} u)(t) = 0$ when $K_1(t) = c_1$. Consequently, the switching times for $(R^{c_0 c_1} u)(t)$ are completely determined by the signal $K(t)$. ///

Next, we address the question of measurability of the relay operator $R^{c_0 c_1}$ as a function of $c_0$, $c_1$. This question is of interest because the complete hysteresis operator will be defined as an integral (continuous superposition of an infinite number of relays) over the variables $(c_0, c_1)$, and measurability is needed to make the integral meaningful. We have the following:

Theorem 2.3. Let $u(t)$ be continuous and such that $K_0(t)$, $K_1(t)$ are piecewise monotone. We assume that all admissible values of $(c_0, c_1)$ make up an open set $\Omega$ in $IR^2$ and $\gamma_0(c_0)$ and $\gamma_1(c_1)$ are continuous as functions of $(c_0, c_1)$ with the metric

$$\rho(\gamma_0(c_0), \gamma_0(c_0')) = \inf_{\substack{x \in \gamma_0(c_0) \\ x' \in \gamma_0(c_0')}} |x - x'|;$$ the metric is defined for $\gamma_1(c_1)$ analogously.

If $i_0 : C \to \{0, 1\}$ is measurable, then for every $t \in [0, T]$, the function $r : C \to \{0, 1\}$ is defined as $r(c_0, c_1) = (R^{c_0 c_1} u)(t)$ is measurable as a function of $(c_0, c_1)$.

<u>Proof.</u> As $K_0(t)$ increases for every $t \in [t_0, t_1)$, the relays that will be turned from "*0*" to "*1*" at time $t$ are contained in the set of relays that have $c_0 \leq K_0(t)$ and those relays have $(c_0, c_1)$ in an open set $\Omega$ ( by assumption), so the set of relays that will be turned from "*0*" to "*1*" is a subset of a measurable set. The set of relays that were in position "*0*" at time $t_0$ has parameters $(c_0, c_1)$ in a measurable set, by the assumption of measurability of $i_0$. Thus the relays that will be turned from "*0*" to "*1*" comprise the intersection of the two sets, the $(c_0, c_1)$ that satisfies $c_0 \leq K_0(t)$, and the values of $(c_0, c_1)$ for which $i_0(c_0, c_1) = 0$. Suppose $K_0(t)$, $K_1(t)$ are piecewise monotone. For each $t \in [t_0, t_1)$, take a partition $t_0 = \tau_0 < \tau_1 < \tau_2 < \cdots < \tau_n = t_1$.
Each $(\tau_i, \tau_{i+1})$ is an interval of monotonicity of $K_0(t)$.
(We can repeat the same argument for $K_1(t)$ on each interval $(\tau_i, \tau_{i+1})$.)
Thus we have measurability.
With piecewise monotonicity of $K_0(t)$ and $K_1(t)$, measurability of $i_0$, and the admissible $(c_0, c_1)$ make up an open set, the outputs will be measurable as functions of $(c_0, c_1)$ for every $t \in [0, T]$. ///

<u>Theorem 2.4.</u> Under the conditions of Theorem 2.3, the function $(R^{c_0 c_1} u)(t)$ has a jointly measurable modification in $(c_0, c_1, t)$.

<u>Proof.</u> If $K_0(\cdot)$ increases from $t_0$ to $t$, then the values of $(c_0, c_1)$ that satisfy $|(R^{c_0 c_1} u)(t) - (R^{c_0 c_1} u)(t_0)| = 1$ consist of the union of two sets:
One is $A_0(t) = \{(c_0, c_1) : K_0(t_0) \leq c_0 \leq K_0(t)\}$: These are the values of $(c_0, c_1)$ that have $(R^{c_0 c_1} u)(t_0) = 0$ and $(R^{c_0 c_1} u)(t) = 1$.
There is an analogous set for those $(c_0, c_1)$ for which $(R^{c_0 c_1} u)(t_0) = 1$ and $(R^{c_0 c_1} u)(t) = 0$.
For every Lebesgue-measurable set $A_0(t) = \{(c_0, c_1) : K_0(t_0) \leq c_0 \leq K_0(t)\}$, define the probability measure

$$P(A_0) = \iint_{A_0} f(x,y)dxdy$$

where $(x,y) \in C$ is all admissible values of $(c_0, c_1)$. $f$ is continuous, $f(x,y) > 0$ for all $(x,y) \in IR^2$, and $\iint_{IR^2} f(x,y)dxdy = 1$.

Set
$$\overline{A} = \{(c_0, c_1) \in C : c_0 = K_0(t_0)\},$$
then
$$A_0 \to \overline{A} = \{(K_0(t_0), c_1) : (K_0(t_0), c_1) \in C\} \text{ as } t \to t_0.$$

The set $\overline{A}$ has two dimensional Lebesgue measure 0, thus also $P(\overline{A}) = 0$. According to the standard results in the theory of stochastic processes, a sufficient condition for existence of a measurable modification of $(R^{c_0 c_1} u)(t)$, as a function of $(c_0, c_1, t)$, is that

$$P\{(c_0, c_1) : |(R^{c_0 c_1} u)(t) - (R^{c_0 c_1} u)(t_0)| = 1\} \to 0 \text{ as } t \to t_0.$$

By continuity of the integral as a function of the domain of integration, $P(A_0) \to 0$ as $t \to t_0$.
Thus there exists a measurable modification. ///

## 3. Bounded variation properties of the hysteresis operator

The properties of bounded variation for the output of the standard Preisach operator have been proved in [BV, V]. In this section, we shall extend the bounded variation properties to our model of two-dimensional hysteresis.

For the function $(R^{c_0 c_1} u):[0,T] \to \{0, 1\}$ and any partition $P := \{0 = t_0 < t_1 < t_2 < \cdots < t_n = T\}$, if the total variation is finite, that is

$$\sup_P \sum_{i=0}^{n-1} \left| (R^{c_0 c_1} u)(t_{i+1}) - (R^{c_0 c_1} u)(t_i) \right| < \infty,$$

$(R^{c_0 c_1} u)(t)$ is of <u>bounded variation</u> over $[0,T]$. And we denote

$$V_0^T (R^{c_0 c_1} u) := \sup_P \sum_{i=0}^{n-1} \left| (R^{c_0 c_1} u)(t_{i+1}) - (R^{c_0 c_1} u)(t_i) \right|$$

(3.1)

<u>Definition 3.1.</u> The function $\omega(u,\delta)$ is said to be the modulus of continuity of a continuous function $u(\cdot)$ for $t_1$, $t_2$ in $[0,T]$ if

$$\omega(u;\delta) := \sup\{r \in IR : |t_1 - t_2| < r \Rightarrow |u(t_1) - u(t_2)| < \delta\}, \text{ for } \delta > 0. \; ///$$

(Note that what we defined as modulus of continuity above is roughly the inverse function of what is defined in some real analysis texts as "modulus of continuity".)

The relay operator $(R^{c_0 c_1} u)(t)$ has following properties:

<u>Theorem 3.1</u> Let $u(\cdot)$ be a continuous piecewise monotone signal on interval $[0,T]$. For an admissible pair of parameter $(c_0, c_1)$, the relay operator $(R^{c_0 c_1} u)(t)$ is piecewise constant and of bounded variation on $[0,T]$. And

$$V_0^T (R^{c_0 c_1} u) \le \frac{T}{\omega(u,\delta)} + 1$$

(3.2)

where $\delta \equiv \inf\{|x - y|: x \in (\gamma_0), y \in (\gamma_1)\} > 0$ and the states of the output of the non-ideal relay are "0" and "1".

Proof. The states of the output of the non-ideal relay $\left(R^{c_0 c_1} u\right)(t)$ are "*0*" and "*1*".
A contribution of 1 to the total variation $V_0^T\left(R^{c_0 c_1} u\right)$ takes place when $u(\cdot)$ changes either from $x \in (\gamma_0)$ to $y \in (\gamma_1)$ or from $y \in (\gamma_1)$ to $x \in (\gamma_0)$. If a change in $u(\cdot)$ from $x \in (\gamma_0)$ to $y \in (\gamma_1)$ or from $y \in (\gamma_1)$ to $x \in (\gamma_0)$ occurs at times $t_1$, $t_2$ in $[0,T]$, then
$$|u(t_1) - u(t_2)| < \delta.$$
So we have
$|u(t_1) - u(t_2)| < \delta$ if $|t_1 - t_2| < r$ for some $r \leq \omega(u;\delta)$; thus,
if $|u(t_1) - u(t_2)| = \delta$, we must have $|t_1 - t_2| \geq \omega(u;\delta)$.
The number of such changes in the values of $u(\cdot)$ from $x \in (\gamma_0)$ to $y \in (\gamma_1)$ or from $y \in (\gamma_1)$ to $x \in (\gamma_0)$ is at most $\dfrac{T}{\omega(u;\delta)}$, and every intermediate variation contributes 1 to the total variation $V_0^T\left(R^{c_0 c_1} u\right)$. Therefore, we have

$$V_0^T\left(R^{c_0 c_1} u\right) \leq \frac{T}{\omega(u,\delta)} + 1. \quad ///$$

Theorem 3.2. If $u(\cdot)$ is continuous and of bounded variation over $[0,T]$, then

$$V_0^T\left(R^{c_0 c_1} u\right) \leq \frac{V_0^T(u)}{\delta}$$

(3.3)

Proof. Suppose $u(\cdot)$ is of bounded variation. Since if changes of $u(\cdot)$ from $x \in (\gamma_0)$ to $y \in (\gamma_1)$ or from $y \in (\gamma_1)$ to $x \in (\gamma_0)$ occur at times $t_0$, $t_1$, $t_2$, ..., $t_n$, we have $n$ jumps between and $x \in (\gamma_0)$ and $y \in (\gamma_1)$. That is, there are $n$ changes of $u(t)$ from $x \in (\gamma_0)$ to $y \in (\gamma_1)$ or from $y \in (\gamma_1)$ to $x \in (\gamma_0)$. Thus, the total number of changes of $u(\cdot)$ from $x \in (\gamma_0)$ to $y \in (\gamma_1)$ or from $y \in (\gamma_1)$ to $x \in (\gamma_0)$ is at most $\dfrac{V_0^T(u)}{\delta}$.
Hence

$$\sum_{i=0}^{n-1} |u(t_{i+1}) - u(t_i)| = n\delta \leq V_0^T(u),$$

$$n \leq \frac{V_0^T(u)}{\delta},$$

Therefore,

$$V_0^T\left(R^{c_0 c_1} u\right) \leq \frac{V_0^T(u)}{\delta}. \quad ///$$

In the next theorem, we state and prove a property of minimum total variation for our model of two-dimensional hysteresis. In the case of the standard Preisach model, this property has been stated in [V].

Theorem 3.3. Let $u(\cdot)$ be an input signal in $C([0,T])$, $\xi \in \{0, 1\}$, and $D_0, D_1$ be two open subsets defined in Definition 2.1. And let $(R(u;\xi))(t)$ be any function that satisfies:

$$(R(u;\xi))(0) = \xi \quad \text{if } u(0) \in D_0 \cap D_1;$$

$$(R(u;\xi))(t) = 0 \quad \text{if } u(t) \notin \overline{D}_1;$$

$$(R(u;\xi))(t) = 1 \quad \text{if } u(t) \notin \overline{D}_0.$$

Let $(R^{c_0 c_1}(u;\xi))(t)$ be the output of the non-ideal relay operator, defined by

$$(R^{c_0 c_1}(u;\xi))(0) = \xi \quad \text{if } u(0) \in D_0 \cap D_1;$$
$$(R^{c_0 c_1}(u;\xi))(t^+) = 1 \quad \text{if } u(t) \in \gamma_0(c_0);$$
$$(R^{c_0 c_1}(u;\xi))(t^-) = 0 \quad \text{if } u(t) \in \gamma_1(c_1);$$
$$(R^{c_0 c_1}(u;\xi))(t^+) = (R^{c_0 c_1}(u;\xi))(t^-) \quad \text{if } u(t) \notin \{\gamma_0(c_0), \gamma_1(c_1)\}.$$

Assume that $(R(u;\xi))(\cdot)$ has finite total variation over $[0,T]$. Then

$$V_0^T(R(u;\xi)) \leq V_0^T(R^{c_0 c_1}(u;\xi))$$

(3.4)

with equality for all $t \in [0,T]$ if and only if $R(u;\xi) \equiv R^{c_0 c_1}(u;\xi)$.

Proof. Since $R(u;\xi)$ takes only the values "0" or "1", a contribution of 1 to the variation $V_0^T(R(u;\xi))$ is made every time $R(u;\xi)$ changes values either from "0" to "1", or from "1" to "0". The output $(R(u;\xi))(t)$ may switch from "0" to "1" or from "1" to "0" only for $u(t) \in D_0 \cap D_1$. The output $(R^{c_0 c_1}(u;\xi))(t)$ may switch either from "0" to "1", or from "1" to "0", only for $u(t) \in \gamma_0(c_0)$ or $u(t) \in \gamma_1(c_1)$. The output $(R(u;\xi))(t)$ may switch from "0" to "1" or from "1" to "0" for values of $u(t) \in D_0 \cap D_1$ as well as for values of $u(t)$ in $\{\gamma_0(c_0), \gamma_1(c_1)\}$. Suppose $(R(u;\xi))(t)$ has switched from "0" to "1" at some time $t_1$ for which $u(t_1) \in D_0 \cap D_1$, and assume also $(R^{c_0 c_1}(u;\xi))(t^-) = 0$; then $(R^{c_0 c_1}(u;\xi))$ will either switch from "0" to "1" at some time $t_2 > t_1$, when $u(t_2) \in \gamma_0(c_0)$ or it will remain at $(R^{c_0 c_1}(u;\xi))(t) = 0$ for $t_1 < t \leq t_3$, where $t_3$ is the first time after $t_1$ that $u(t)$ reaches the value $u(t) \in \gamma_1(c_1)$; in the first case, the total contribution to both $V_0^{t_2}(R(u;\xi))$ and $V_0^{t_2}(R^{c_0c_1}(u;\xi))$ is the same, that is, a contribution of $+1$; in the later case, the contribution to $V_0^{t_3}(R(u;\xi))$ is 2, and the contribution to $V_0^{t_3}(R^{c_0c_1}(u;\xi))$ is 0. For the strict inequality: If $(R^{c_0c_1}(u;\xi))(t_1) = 0$ and $(R(u;\xi))(t_1) = 0$, and if, at time $t_1$, $u(t_1) \in D_0 \cap D_1$ and $(R(u;\xi))$ switches from "*0*" to "*1*" at time $t_1$, then the contribution to $V_0^{t_1}(R(u;\xi))$ is 1, whereas the contribution to $V_0^{t_1}(R^{c_0c_1}(u;\xi))$ is 0. So, at time $t_1$, $V_0^{t_1}(R(u;\xi)) \le V_0^{t_1}(R^{c_0c_1}(u;\xi))$. By considering all other cases in a similar way, we can find that, if $(R(u;\xi))(t) \ne (R^{c_0c_1}(u;\xi))(t)$ then there exists some $t$ such that $V_0^t(R(u;\xi)) \le V_0^t(R^{c_0c_1}(u;\xi))$. ///

## 4. The wiping-out property and related issues in two-dimensional hysteresis

We denote by $C$ the set of all admissible values of $(c_0, c_1)$ (cf. Definition 2.2). We assume that C is a closed bounded set in $IR^2$ and is equal to the closure of its own interior. The corresponding admissible values of $(K_0(t), K_1(t))$ form another set in $IR^2$. According to case (ii) of Theorem 2.2, for admissible values of $(K_0(t), K_1(t))$, it is not possible to find $(c_0, c_1) \in C$ such that both $c_0 < K_0(t)$ and $c_1 < K_1(t)$.

Theorem 4.1. Under the properties (1)-(4) in section 2, the hysteresis operator

$$(Hu)(t) = \iint_{(c_0, c_1)} w(c_0, c_1)(R^{c_0, c_1} u)(t) dc_0 dc_1$$

where $w(c_0, c_1)$ is weight function

(4.1)

depends only on the history of the local extrema of the two signals $K_0(t)$ and $K_1(t)$, and the history of $K_0(\cdot)$ and $K_1(\cdot)$ after the most recent local extrema.

Proof. At a point $t$ that is not a point of local extremum of either $K_0(t)$ or $K_1(t)$, we have four possibilities for the pair $(K_0(t), K_1(t))$:

(i) $K_0(\cdot)$ increases at $t$, and $K_1(\cdot)$ also increases at $t$.
(ii) $K_0(\cdot)$ increases at $t$, and $K_1(\cdot)$ decreases at $t$.
(iii) $K_0(\cdot)$ decreases at $t$, and $K_1(\cdot)$ increases at $t$.
(iv) $K_0(\cdot)$ decreases at $t$, and $K_1(\cdot)$ also decreases at $t$.

First, we examine the behavior of $(Hu)(t)$ in case (i).
Let $\bar{t}_0$ be the most recent local extremum $K_0(\cdot)$ before time of $t$, and $\bar{t}_1$ be the most recent local extremum of $K_1(\cdot)$ before time $t$. At the time $\bar{t}_0$ certain relays were in the "*0*" position, and certain relays were in the "*1*" position. For $\bar{t}_0 < t' \le t$, since we are in the case of increasing $K_0(\cdot)$, all the relays that were in the "*0*" position at time $\bar{t}_0$ and have $c_0 = K_0(t')$, will be turned to the "*1*" position. Similarly, for $\bar{t}_1 < t'' \le t$, since $K_1(\cdot)$ increases, those relays that were in the "*1*" position at time $\bar{t}_1$ and have $c_1 = K_1(t'')$ will be turned to the "*0*" position. In this way, the output states of all relays at time depend on the output states at times $\bar{t}_0$ and $\bar{t}_1$, the values of $K_0(t)$ and $K_1(t)$, and the information that $K_0(\cdot)$ is increasing over $(\bar{t}_0, t]$ and $K_1(\cdot)$ is also increasing over $(\bar{t}_1, t]$.

The behavior of the relays, and thus also the behavior of $(Hu)(t)$ in cases (ii) and (iii), and (iv), can be examined in a similar way. ///

Remark 4.1. We have established that a general signal $u(t) = (u_1(t), u_2(t))$ leads to a signal $K(t) = (K_0(t), K_1(t))$ that suffices to determine the hysteresis operator. The question arises, whether a signal $K(t)$ can be used to determine the original input signal $u(t)$. In general, this is not possible, as it can be shown by simple examples. If, however, we restrict attention to a domain $\Omega$ with the property that, within $\Omega$, a curve of the group $(\gamma_0)$ and a curve of the group $(\gamma_1)$ can intersect at at most one point, then, for $u(t)$ remaining in $\Omega$, the correspondence between $u(t)$ and $K(t)$ is one-to-one. ///

In the case of the standard Preisach model, the wiping-out property is formulated in terms of the dominant local extrema of the input signal $u(\cdot)$; cf. [M1].
For our two-dimensional model, the analogue of dominant local extrema is the concept of dominant reversal points. We proceed to formulate the relevant definitions.

Definition 4.1. A reversal point of the two-dimensional signal $u(\cdot)$ is a point $u(t_1)$ where $t_1$ is a time-instant of a local maximum of either $K_0(\cdot)$ or $K_1(\cdot)$.
A reversal point will be called <u>reversal point of type 0</u> if it corresponds to a local maximum of $K_0(\cdot)$, and a <u>reversal point of type 1</u> if it corresponds to a local maximum of $K_1(\cdot)$. ///

Definition 4.2. A signal $u(\cdot)$ will be called <u>regular</u> if $u(\cdot)$ is continuous and the two signals $K_0(\cdot)$ and $K_1(\cdot)$ are both piecewise monotone. ///

Definition 4.3. The dominant reversal points of a regular signal $u(\cdot)$ are defined as follows: the first local maximum of either $K_0(\cdot)$ or $K_1(\cdot)$ is designated as a dominant reversal point; inductively, if $t_i$ is a dominant reversal point and a local maximum of $K_j(\cdot)$ $(j = 0, 1)$, then the next dominant reversal point $t_{i+1}$, $i = 1, 2, 3, \cdots$, is a point of local maximum of $K_{j+1}(\cdot)$ (addition of subscripts is modulo 2) with the property that $K_j(t) \leq K_j(t_i)$ for all $t \in [t_i, t_{i+1}]$. ///

The wiping – out property for the two dimensional case of hysteresis operator can be stated in the following form:

Theorem 4.2. The output of the two-dimensional hysteresis operator (4.1) at time $t$ for regular input signals $u(\cdot)$, depends only on the history of dominant reversal points of $u(\cdot)$ prior to time $t$, and the history of $u(\cdot)$ after the last dominant reversal point until time $t$.

Proof. We prove this by induction on the number $N$ of dominant reversal points of $u(\cdot)$ in the interval $[0,T]$.

If $N=1$, and there is only one dominant reversal point $t_1$, which may assume, without loss of generality, to be of type *0*, then this means that $t_1$ is a point of local maximum of $K_0$, and also a point of global maximum of $K_0$, and there is no local maximum of the signal $K_1$ in $[0,T]$. Then, by time $t_1$, all the relays $R^{c_0 c_1}$ with $c_0 \leq K_0(t_1)$, have been turned to the output state "*1*". In this case, the conclusion of the theorem holds, since we know the state of all relays at time $t_1$, and the hysteresis output at time $t$ depends only on the history of the input signal after time $t_1$.

Inductively, if, by time $t_k$, where $t_k$ is a dominant reversal point, the hysteresis output depends only on the dominant reversal points up to and including time $t_k$, then we examine it at the next dominant reversal point $t_{k+1}$.

Without loss of generality, we assume that $t_k$ is a dominant reversal point of type *0*, so that, by definition, $t_{k+1}$ is a dominant reversal point of type *1*. At time $t_{k+1}$, all the relays $R^{c_0 c_1}$ with $c_1 \leq K_1(t_{k+1})$ will be turned to the output state "*0*". For $t \in [t_k, t_{k+1}]$, we have $K_0(t) \leq K_0(t_k)$, thus there will be no relays $R^{c_0 c_1}$ that have been turned from output state "*0*" to output state "*1*" for time $t \in [t_k, t_{k+1}]$. Thus, the relays that are turned off from output state "*1*" to output state "*0*" over $t \in [t_k, t_{k+1}]$ are those of the relays $R^{c_0 c_1}$ that were in output state "*1*" at time $t_k$ and satisfy $c_1 \leq K_1(t_{k+1})$. Therefore, the hysteresis output at time $t_{k+1}$ depends only on the state of relays up to time $t_k$ and the value of the signal $u(\cdot)$ between $t_k$ and $t_{k+1}$.

The induction is completed, and the conclusion holds for any number $N$ of dominant reversal points in $[0,T]$. ///

5. Identification methods

The identification of the weight function $w(\alpha, \beta)$ in a standard Preisach model is one of the basic problems of the theory of hysteresis. This problem has been solved in [M1], and the solution relies on information about the first-order transition curves.

In the case of the two-dimensional hysteresis models introduced in the present paper, the identification problems become more complicated. There are two relevant identification problems that arise naturally in the context of the two-dimensional hysteresis operators introduced in this paper.

(i) If the families of curves, $\gamma_0(c_0)$ and $\gamma_1(c_1)$, introduced in section 2, are completely known, we want to identify the weight function $w(c_0, c_1)$ of the two-dimensional hysteresis operator $H$, by using observations of the output $H(t)$ for suitably chosen two-dimensional input signals $u(\cdot)$.

(ii) The families of curves $\gamma_0(c_0)$ and $\gamma_1(c_1)$ may not be known. In that case, the identification problem consists of two parts: identification of the curves $\gamma_0(c_0)$, $\gamma_1(c_1)$, on the basis observations of the output $(Hu)(t)$ for suitably chosen two-dimensional input signals $u(t)$, and identification of the weight function $w(c_0, c_1)$.

For the first identification problem, we postulate the following properties.

(P1) $\Omega$ has nonempty boundary $\partial \Omega$.

(P2) There exists a part $\Gamma$ of $\partial \Omega$ with the property that, if an input signal $u(\cdot)$ starts outside $\Omega$ and enters $\Omega$ through $\Gamma$, then all relay outputs $\left(R^{c_0 c_1} u\right)(t)$ are $0$ when

$$t = t_e \equiv \inf\{t > t_0 : u(t) \notin I\!R^2 - \overline{\Omega}\}$$

(that is, $t_e$ is the time of entrance of $u(\cdot)$ into $\Omega$ through the part $\Gamma$ of $\partial \Omega$).

(P3) There exists a simple $C^1$ curve $(k)$ that intersects $\Gamma$ and also intersects all curves $\gamma_0(c_0)$ and $\gamma_1(c_1)$ transversally.

(P4) For every point $P \in (k)$, if $c_0(P) < c_0(P')$, then the pair $\left(c_0(P'), c_1(P)\right)$ is admissible.

(P5) All points $P$, $Q$ on $(k)$ can be parametrized by arc length $s$ on $(k)$ in such a way that

$$P < Q \Leftrightarrow s(P) < s(Q).$$

(P6) Let $s_0(\tilde{c}_0)$ be the value of $s$ for which $c_0(P(s)) = \tilde{c}_0$ and $s_1(\tilde{c}_1)$ be the values of $s$ for which $c_1(Q(s)) = \tilde{c}_1$. (where $c_0(P(s))$ and $c_1(Q(s))$ are $C^1$ as functions of $s$.) Then, for $s_0 > s_1$, the Jacobian

$$J(s_0, s_1) = \left| \frac{\partial(c_0(P(s_0)), c_1(Q(s_1)))}{\partial(s_0, s_1)} \right|$$

is never zero.

Properties (P1), (P2), and (P3) are illustrated in Figure 5.1.
We have:

<u>Lemma 5.1.</u> If $\overline{c_0} = c_0(P)$ and $c_1' = c_1(P')$ are such that the pair $(\overline{c_0}, c_1')$ is admissible, with distinct points $P$, $P'$ on $(k)$, then we have

$$c_0(P') < c_0(P).$$

<u>Proof.</u> If $(\overline{c_0}, c_1')$ is an admissible pair, then $D_0(\overline{c_0}) \cup D_1(c_1') = \Omega$.
Therefore, $\gamma_0(\overline{c_0}) \subseteq D_1(c_1')$, since $D_0(\overline{c_0})$ is an open set and $\gamma_0(\overline{c_0}) \subseteq \partial D_0(\overline{c_0})$, thus $\gamma_0(\overline{c_0}) \cap D_0(\overline{c_0}) = \phi$, and consequently,
$$\gamma_0(\overline{c_0}) \subseteq \Omega \setminus D_0(\overline{c_0}) \subseteq D_1(c_1').$$
Similarly, $\gamma_1(c_1') \subseteq D_0(\overline{c_0})$.
Let $c_0' = c_0(P')$, $\overline{c_1} = c_1(P)$.
From condition (4) of section 2, we have:
   if $D_0(\overline{c_0}) \subsetneq D_0(c_0')$, then $D_1(c_1') \subsetneq D_1(\overline{c_1})$.
We must have $c_0' \neq \overline{c_0}$, since, by our assumption, the curve $(k)$ intersects each curve $\gamma_0(c_0)$ at exactly one point, and therefore it is impossible to have $c_0' = \overline{c_0}$ since that would imply that the two distinct points $P$, $P'$ on $(k)$ lie on the same curve $\gamma_0(c_0)$.
Then, by property (2) of section 2, exactly one of the two statements must be true:

$$D_0(\overline{c_0}) \subsetneq D_0(c_0'), \text{ or } D_0(c_0') \subsetneq D_0(\overline{c_0}).$$

If $D_0(\overline{c_0}) \subsetneq D_0(c_0')$, and $(\overline{c_0}, c_1')$ is admissible, that is,
   $D_0(\overline{c_0}) \cup D_1(c_1') = \Omega$, then also
   $D_0(c_0') \cup D_1(c_1') = \Omega$,

thus $(c_0',c_1')$ is admissible. But $P'$ lies on both $\gamma_0(c_0')$ and $\gamma_1(c_1')$, and therefore $(c_0',c_1')$ can not be admissible since the curves $\gamma_0$, $\gamma_1$ corresponding to an admissible pair can not intersect.

Therefore, $D_0(c_0') \subsetneq D_0(\overline{c_0})$, which by condition (3) of section 2, is equivalent to

$$c_0(P') < c_0(P). \;///$$

From Lemma 5.1 and property (P4), we have:

<u>Corollary 5.1.</u> For any two distinct points $P$, $P'$ on $(k)$, the pair $\big(c_0(P'),c_1(P)\big)$ will be admissible if and only if $c_0(P) < c_0(P')$. ///

<u>Definition 5.1</u> For any two distinct points $P$, $P'$ on $(k)$, we define the total order $\leq$ by

$$P < P' \text{ if and only if } c_0(P) < c_0(P'), \text{ and}$$
$$P \leq P' \text{ if and only if either } P < P' \text{ or } P = P'. \;///$$

We have:

<u>Lemma 5.2.</u> If $u(\cdot)$ is a continuous input signal restricted to take values on the curve $(k)$, then the relay $(R^{\overline{c_0 c_1}} u)(\cdot)$ is equivalent to a non-ideal relay of the standard Preisach type, in the following sense: for $P$, $Q$ on $(k)$, with $P > Q$, define a relay $(S^{PQ} u)(\cdot)$ by

$$(S^{PQ} u)(t^+) = (S^{PQ} u)(t^-), \text{ if } Q < u(t) < P;$$
$$(S^{PQ} u)(t^+) = 1, \text{ if } u(t) \geq P;$$
$$(S^{PQ} u)(t^+) = 0, \text{ if } u(t) \leq Q.$$

<u>Proof.</u> We show that if $\overline{c_0} = c_0(P)$ and $c_1' = c_1(Q)$, then $(R^{\overline{c_0 c_1'}} u)(t) = (S^{PQ} u)(t)$ at points $t$, where $(R^{\overline{c_0 c_1'}} u)(t)$ is continuous, and
$(R^{\overline{c_0 c_1'}} u)(t^+) = (S^{PQ} u)(t^+)$ at point $t$ of jump discontinuity of $(R^{\overline{c_0 c_1'}} u)$.
We have, by Corollary 5.1, that there is a one-to-one correspondence between pairs $(P,Q)$ on $(k)$ with $P > Q$ and admissible pairs $(\overline{c_0}, c_1')$, where $\overline{c_0} = c_0(P)$ and $c_1' = c_1(Q)$. We have:

$$u(t) \in \gamma_0(\overline{c_0}) \Leftrightarrow u(t) = P, \text{ and } u(t) \in \gamma_1(c_1') \Leftrightarrow u(t) = Q,$$

since by assumption $(k)$ intersects each curve of the family $(\gamma_0)$ at exactly one point, and also intersects each curve of the family $(\gamma_1)$ at exactly one point.

The condition $u(0) \in D_0(\overline{c_0}) \cap D_1(c_1')$ is equivalent to $Q < u(0) < P$; this is true because
$$u(0) > Q \Leftrightarrow c_0(Q) < c_0(u(0))$$
and by property (4) of section 2, this is equivalent to
$$c_1(Q) > c_1(u(0)),$$
that is,
$$c_1(u(0)) < c_1'.$$
Since $u(0) > Q$, the pair $(c_0(u(0)), c_1')$ is admissible, thus, by the same argument as in the proof of Lemma 5.1,
$$\gamma_0(c_0(u(0))) \subseteq D_1(c_1'),$$
thus
$$u(0) \in D_1(c_1').$$
Similarly, $u(0) < P$ implies $u(0) \in D_0(\overline{c_0})$. Thus the double inequality $Q < u(0) < P$ implies $u(0) \in D_0(\overline{c_0}) \cap D_1(c_1')$. By the definition of the total order on $P$ and condition (4) of section 2,
$$\text{if } u(0) \in D_0(\overline{c_0}) \cap D_1(c_1'), \text{ then } Q < u(0) < P.$$
Thus the condition $u(0) \in D_0(\overline{c_0}) \cap D_1(c_1')$ is equivalent to $Q < u(0) < P$.

The condition $u(t) \in \gamma_0(\overline{c_0})$ is equivalent to $u(t) = P$, since $(k)$ intersects $\gamma_0(\overline{c_0})$ at the single point $P$.

Similarly, the condition $u(t) \in \gamma_1(c_1')$ is equivalent to $u(t) = Q$.

Thus the definition of $(S^{PQ}u)(\cdot)$ given in this lemma is equivalent to the characterization of $(R^{\overline{c_0}c_1'}u)(\cdot)$ stated in the Theorem 3.3.

Consequently,
$$(R^{\overline{c_0}c_1'}u)(\cdot) = (S^{PQ}u)(\cdot). \quad ///$$

Definition 5.2. The space $I_2((k))$ consists of two-dimensional continuous input signals that take values on the curve $(k)$. ///

A consequence of Lemma 5.2 is the following:

Theorem 5.1. If we set $W(P,Q) = w(c_0(P), c_1(Q))$ for $P$, $Q$ on $(k)$ with $P > Q$, and if $u(\cdot)$ is an input signal in $I_2((k))$, then the two-dimensional hysteresis operator $(Hu)(t)$ is completely characterized by the function $W$.

Proof. Since $(k)$ intersects each $\gamma_0$ and $\gamma_1$ at exactly one point, for all $c_0$ and $c_1$, we define that $P(c_0)$ is the point of intersection of $(k)$ with $\gamma_0(c_0)$, and $Q(c_1)$ is the point of intersection of $(k)$ with $\gamma_1(c_1)$. Then

$$w(c_0,c_1) = W(P(c_0),Q(c_1)) \text{ and}$$
$$(c_0,c_1) \text{ is an admissible} \Leftrightarrow P(c_0) > Q(c_1)$$

For $u \in I_2((k))$,

$$(Hu)(t) = \iint_{(c_0,c_1)\text{admissible}} w(c_0,c_1)(R^{c_0 c_1}u)(t)dc_0 dc_1$$

$$= \iint_{(c_0,c_1)\text{admissible}} W(P(c_0),Q(c_1))(S^{P(c_0)Q(c_1)}u)(t)dc_0 dc_1$$

thus $(Hu)(t)$ depends only on the weight function $W(P,Q)$. ///

Now, we invoke properties (P5) and (P6). The integral

$$\iint_{(c_0,c_1)\text{admissible}} W(P(c_0),Q(c_1))(S^{P(c_0)Q(c_1)}u)(t)dc_0 dc_1$$

is transformed to an integral in variables $(s_0,s_1)$, according to

$$\iint_{(c_0,c_1)\text{admissible}} W(P(c_0),Q(c_1))(S^{P(c_0)Q(c_1)}u)(t)dc_0 dc_1$$

$$= \iint_{s_0 > s_1} W(P(s_0),Q(s_1))J(s_0,s_1)(S^{P(s_0)Q(s_1)}u)(t)ds_0 ds_1$$

where $J(s_0,s_1)$ is the Jacobian defined in property (P6).

The expression

$$\Phi(s_0,s_1) \equiv W(P(s_0),Q(s_1))J(s_0,s_1)$$

(5.1)

can be identified by the same methods as in the standard Preisach model in [M1]. For completeness, we outline that method. First, an input signal $u(\cdot)$ is steadily increased until it achieves the value $s_0$, then it is decreased until it achieves the value $s_1$.

The difference
$$\psi(s_0, s_1) := (Hu)(s_0) - (Hu)(s_1)$$
is the first-order transition function.
Then $\Phi(s_0, s_1)$ is given by

$$\Phi(s_0, s_1) = -\frac{\partial^2 \psi(s_0, s_1)}{\partial s_0 \partial s_1}$$

Then we have:

Theorem 5.2. The function $W(P, Q)$ can be identified by using

$$W(\overline{P}(s_0), \overline{Q}(s_1)) = -\frac{1}{J(s_0, s_1)} \frac{\partial^2 \psi(s_0, s_1)}{\partial s_0 \partial s_1}$$

The original continuous weight function $w(c_0, c_1)$ can be identified by using

$$w(c_0, c_1) = W(P(c_0), Q(c_1)).$$

Proof. As outlined above, the function $\Phi(s_0, s_1)$ satisfies

$$\Phi(s_0, s_1) = -\frac{\partial^2 \psi(s_0, s_1)}{\partial s_0 \partial s_1}.$$

By property (P6), the Jacobian $J(s_0, s_1)$ is never zero for $s_0 > s_1$. Thus,

$$W(\overline{P}(s_0), \overline{Q}(s_1)) = -\frac{1}{J(s_0, s_1)} \frac{\partial^2 \psi(s_0, s_1)}{\partial s_0 \partial s_1}$$

The equality
$$w(c_0, c_1) = W(P(c_0), Q(c_1))$$
for admissible $(c_0, c_1)$ is the definition of $W(P(c_0), Q(c_1))$. ///

For the second identification problem, namely the problem of simultaneous identification of the weight function $w(c_0, c_1)$ and the two families of curves $(\gamma_0(c_0), \gamma_1(c_1))$, we assume, in addition to previous properties, that there exists a family $\{(k_\xi) : \xi \in \Xi\}$ of curves, each curve $(k_\xi)$ satisfying properties (P1)-(P6), such that the curves $(k_\xi)$ do not intersect each other, the sets

$$\{(k_\xi) \cap \gamma_i(c_i) : \xi \in \Xi\}, \quad i = 0, 1$$

are dense in $\gamma_i(c_i)$, for each value of $c_i$, and each set

$$\Gamma_{i,\xi}(c_i) \equiv (k_\xi) \cap \gamma_i(c_i)$$

can be expressed as

$$s = G_i(\xi, c_i), \text{ for } i \text{ taking values in } \{0,1\},$$

where $s$ represents arc length on $(k_\xi)$. These properties are illustrated in Figure 5.2.

Let $\Phi_\xi(s_0, s_1)$ be defined, for $s_0 > s_1$, as the function $\Phi(s_0, s_1)$ of (5.1) defined for the hysteresis operator applied to input signals restricted to lie on $(k_\xi)$, and let $J_\xi(s_0, s_1)$ be the corresponding Jacobian for each curve $(k_\xi)$. Let $W_\xi(P,Q)$ be the corresponding function $W(P,Q)$ defined on $(k_\xi)$.

In this problem, we have a choice for representing the parameters $c_0$, $c_1$ as functions of $s_0$, $s_1$ on each $(k_\xi)$. The simplest choice is $c_0 \equiv s_0$, $c_1 \equiv s_1$.

With this parametrization, we have $J_\xi(s_0, s_1) \equiv 1$ for all $s_0$, $s_1$ with $s_0 > s_1$, and consequently,

$$w(c_0, c_1) = w(s_0, s_1)$$

$$= W_\xi(\overline{P_\xi(s_0)}, \overline{Q_\xi(s_1)})$$

$$= W_\xi(P_\xi(c_0), Q_\xi(c_1))$$

$$= \Phi_\xi(s_0, s_1).$$

We have:

<u>Theorem 5.3.</u> Under the conditions stated above, the curves $\gamma_0(c_0) \equiv \gamma_0(s_0)$, $\gamma_1(c_1) \equiv \gamma_1(s_1)$ are described, respectively, in implicit form by the equations

$$\Phi_\xi(s, s_1) = w(s_0, s_1) \text{ for arbitrary but fixed } s_1,$$
$$\Phi_\xi(s_0, s) = w(s_0, s_1) \text{ for arbitrary but fixed } s_0.$$

Proof. $\gamma_0(s_0)$ is described by an equation

$$s = G_0(\xi, s_0).$$

The intersection of $\gamma_0(s_0)$ with $(k_\xi)$ is the point $\overline{P_\xi}(s_0)$ on $(k_\xi)$. Thus, for $s = G_0(\xi, s_0)$, we must have

$$\Phi_\xi(s, s_1) = w(s_0, s_1)$$

which is an implicit equation that describes

$$\{\gamma_0(s_0) \cap (k_\xi) : \xi \in \Xi\},$$

which is a set that is dense in $\gamma_0(s_0)$ by assumption, thus the equation $\Phi_\xi(s, s_1) = w(s_0, s_1)$ describes the curve $\gamma_0(c_0)$ as an implicit equation. Similarly, for arbitrary but fixed $s_0$, the equation

$$\Phi_\xi(s_0, s) = w(s_0, s_1)$$

with arbitrary but fixed $s_1$, describes the set

$$\{\gamma_1(s_1) \cap (k_\xi) : \xi \in \Xi\},$$

and thus, by the density condition, the curve $\gamma_1(s_1)$, in the form of an implicit equation in the variables $s$, $\xi$. ///